\documentstyle [12pt]{article}

\def\qmod#1#2{{\hbox{}^{\displaystyle{#1}}}\!\big/\!\hbox{}_{
\displaystyle{#2}}}

\newfam\msbfam
\font\twelmsb=msbm10 at 12pt
\font\tenmsb=msbm10
\font\sevenmsb=msbm10 at 7pt
\font\fivemsb=msbm10 at 5pt

\textfont\msbfam=\twelmsb
\scriptfont\msbfam=\sevenmsb
\scriptscriptfont\msbfam=\fivemsb

\def\Bbb{\fam\msbfam\tenmsb}

\def\C{{\Bbb C}}

\def\H{{\Bbb H}}

\def\R{{\Bbb R}}
\def\Z{{\Bbb Z}}

\def\EE{{\cal E}}
\def\FF{{\cal F}}

\def\MM{{\cal M}}

\def\SS{{\cal S}}

\def\union{\mathop{\bigcup}}
\def\qed {\hfill\vrule height6pt width6pt depth0pt \bigskip}

\def\map{\longrightarrow}
\def\textmap#1{\mathop{\vbox{\ialign{
                                ##\crcr
    ${\scriptstyle\hfil\;\;#1\;\;\hfil}$\crcr
    \noalign{\kern-1pt\nointerlineskip}
    \rightarrowfill\crcr}}\;}}

\def\textlmap#1{\mathop{\vbox{\ialign{
                                ##\crcr
    ${\scriptstyle\hfil\;\;#1\;\;\hfil}$\crcr
    \noalign{\kern-1pt\nointerlineskip}
    \leftarrowfill\crcr}}\;}}

\newfam\meuffam
\font\tenmeuf=eufm10
\font\sevenmeuf=eufm7
\font\fivemeuf=eufm5

\textfont\meuffam=\tenmeuf
\scriptfont\meuffam=\sevenmeuf
\scriptscriptfont\meuffam=\fivemeuf

\def\germ{\fam\meuffam\tenmeuf}

\def\c{{\germ c}}

\def\hh{{\germ h}}

\def\p{{\germ p}}

\begin{document}
\def\Pr{{\rm Pr}}
\def\tr{{\rm Tr}}
\def\End{{\rm End}}
\def\Spin{{\rm Spin}}
\def\U{{\rm U}}
\def\SU{{\rm SU}}
\def\SO{{\rm SO}}
\def\PU{{\rm PU}}
\def\spin{{\rm spin}}
\def\u{{\rm u}}
\def\su{{\rm su}}
\def\so{{\rm so}}
\def\pu{{\rm pu}}
\def\Pic{{\rm Pic}}
\def\NS{{\rm NS}}
\def\deg{{\rm deg}}
\def\Hom{{\rm Hom}}
\def\h{{\germ h}}
\def\Herm{{\rm Herm}}
\def\Vol{{\rm Vol}}
\def\pf{{\bf Proof: }}
\def\id{{\rm id}}
\def\i{{\germ i}}
\def\im{{\rm im}}
\def\rk{{\rm rk}}
\def\ad{{\rm ad}}
\def\h{{\bf H}}
\def\coker{{\rm coker}}
\def\dv{\bar\partial}
\def\dva{\bar\partial_A}
\def\da{\partial_A}
\def\p{\partial\bar\partial}
\def\pa{\partial_A\bar\partial_A}
\def\Dr{\hskip 4pt{\not}{D}}
\newtheorem{sz}{Satz}[section]
\newtheorem{th}[sz]{Theorem}
\newtheorem{pr}[sz]{Proposition}
\newtheorem{re}[sz]{Remark}
\newtheorem{co}[sz]{Corollary}
\newtheorem{dt}[sz]{Definition}
\newtheorem{lm}[sz]{Lemma}
\newtheorem{cl}[sz]{Claim}
\begin{titlepage}
\pagestyle{empty}
\title{Quaternionic Monopoles}

\author{Christian Okonek\thanks{Partially supported by: AGE-Algebraic
Geometry in Europe,
contract No ERBCHRXCT940557 (BBW 93.0187), and by  SNF, nr. 21-36111.92}
 \and     Andrei
Teleman$^*$\\ \\ Mathematisches Institut  Universit\"at
Z\"urich\\ Winterthurerstrasse 190, CH-8057 Z\"urich \\  \\ }

 \date{May 25, 1995 }
\maketitle
\begin{abstract}
We present the simplest non-abelian version of Seiberg-Witten theory:
Quaternionic monopoles.
These monopoles  are associated with
$Spin^h(4)$-structures on 4-manifolds and form finite-dimensional moduli
spaces.
On a K\"ahler surface the quaternionic monopole equations decouple and lead
to the projective vortex equation for holomorphic pairs. This vortex
equation comes from a moment map and gives rise to a new complex-geometric
stability concept.
The moduli spaces of quaternionic monopoles on K\"ahler surfaces have two
closed subspaces, both naturally isomorphic with moduli spaces of
canonically stable holomorphic pairs. These components intersect along
Donaldsons instanton space and can be compactified with Seiberg-Witten
moduli spaces. This should provide  a link between the two corresponding
theories.

\end{abstract}

\end{titlepage}

\setcounter{section}{-1}
\section{Introduction}

Recently, Seiberg and Witten [W] introduced new 4-manifold invariants,
essentially by
counting  solutions of the monopole equations. The new invariants have
already found nice
applications, like e.g. in the proof of the Thom conjecture [KM] or in a
short proof of the Van
de Ven conjecture [OT2].  In this paper we introduce and study the simplest
and the most
natural non-abelian version of the Seiberg-Witten monopoles, the
quaternionic monopoles.

Let $(X,g)$ be an oriented Riemannian manifold of dimension 4. The
structure group $SO(4)$
has as natural extension the quaternionic spinor group
$Spin^h(4):=Spin(4)\times_{\Z_2}Sp(1)$:
$$1\map Sp(1)\map Spin^h(4)\map SO(4)\map 1 \ .$$
The projection onto the second factor $Sp(1)=SU(2)$ induces a "determinant map"
$\delta:Spin^h(4)\map PU(2)$.

A $Spin^h(4)$-structure on $(X,g)$ consists of a $Spin^h(4)$-bundle over
$X$  and an
isomorphism of its $Sp(1)$-quotient with the (oriented) orthonormal frame
bundle of
$(X,g)$. Given a $Spin^h(4)$-structure on $X$, one has a one-one
correspondence between
$Spin^h$-connections projecting onto the Levi-Civita connection and
$PU(2)$-connections in
the associated "determinant"  $PU(2)$-bundle. The quaternionic monopole
equations are:
$$\left\{\begin{array}{ccc}\Dr_{A}\Psi&=&0\ \ \ \ \\
\Gamma(F_{A}^+)&=&(\Psi\bar\Psi)_0 \ \ ,
\end{array}\right.\eqno{ }$$
 where $A$ is a $PU(2)$-connection in the "determinant" of the
$Spin^h(4)$-structure  and
$\Dr_A$ the induced Dirac operator; $\Psi$ is a positive quaternionic
half-spinor. The Dirac
operator satisfies the crucial Weitzenb\"ock formula :
$$\Dr_A^2=\nabla_{\hat A}^*\nabla_{\hat A}+\Gamma(F_A)+\frac{s}{4}\id$$
It can be used to show that the solutions of the quaternionic monopole
equations are
the absolute minima of a certain functional, just  like in the
$Spin^c(4)$-case [JPW].

The moduli space of quaternionic monopoles associated with a fixed
$Spin^h(4)$-structure
${\germ h}$ is a real analytic space of virtual dimension
$$m_{\germ h}=-\frac{1}{2}(3p_1+3e+4\sigma)\ .$$
 Here $p_1$ is the first Pontrjagin class of
the determinant, $e$ and $\sigma$ denote the Euler characteristic  and the
signature of $X$.

Note that $m_{\germ h}$ is an even integer iff $X$ admits an almost complex
structure.

The moduli spaces of quaternionic monopoles  contain the Donaldson
instanton moduli spaces as
well as the classical Seiberg-Witten moduli spaces, which suggests that
they could provide a
method of comparing the two theories. We study the analytic structure
around the Donaldson
moduli space.

Much more can be said if the holonomy of $(X,g)$ reduces to $U(2)$, i.e. if
$(X,g)$ is a
K\"ahler surface. In this case we use the canonical $Spin^c(4)$-structure  with
$\Sigma^+=\Lambda^{00}\oplus\Lambda^{02}$ and $\Sigma^-=\Lambda^{01}$ as spinor
bundles. The data of a $Spin^h(4)$-structure ${\hh}$ in $(X,g)$ is then
equivalent to the data
of a Hermitian 2-bundle $E$ with $\det E=\Lambda^{02}$. The determinant
$\delta(\hh)$ coincides with the $PU(2)$-bundle $P(E)$ associated with $E$.
A positive spinor
can be written as
$\Psi=\varphi+\alpha$, where
$\varphi\in A^0(E^{\vee})$ and $\alpha\in A^{02}(E^{\vee})$ are
$E^{\vee}$-valued forms.
To give a $PU(2)$-connection in $P(E)$ means to give a $U(2)$-connection
in $E$ inducing the
Chern connection in $\Lambda^{02}$, or equivalently, a $U(2)$-connection $C$ in
$E^{\vee}$
inducing the Chern connection in $K_X=\Lambda^{20}$.  A pair
$(C,\varphi+\alpha)$ solves the
quaternionic monopole equation iff $C$ is a connection of type $(1,1)$, one
of $\alpha$ or
$\varphi$ vanishes while the other is $\bar\partial_C$-holomorphic, and a
certain projective
vortex equation is satisfied. This shows that in the K\"ahler case the
moduli space decomposes as
a union of two Zariski closed subspaces intersecting along the Donaldson
locus. The two subspaces
are interchanged by a natural real analytic involution, whose fixed point
set is precisely the
Donaldson moduli space.

The projective vortex equation comes from a moment map which corresponds to
a new stability
concept for pairs
$({\cal E},\varphi)$ consisting of a holomorphic bundle ${\cal E}$ with
canonical determinant
$\det{\cal E}={\cal K}_X$ and a holomorphic section $\varphi$. We call such
a pair canonically
stable iff either  ${\cal E}$ is stable, or $\varphi\ne 0$ and the
divisorial component $D_\varphi$
of the zero locus satisfies the inequality
$$c_1\left({\cal O}_X(D_\varphi )^{\otimes2}\otimes K_X^{\vee}\right) \cup
[\omega_g]<0 \ \ .$$

 Our main result identifies the moduli spaces of irreducible quaternionic
monopoles on a
K\"ahler surface with the algebro-geometric moduli space of canonically
stable pairs.

In the algebraic case, moduli spaces of quaternionic monopoles can easily
be  computed using our
main result (Theorem 7.3) and Lemma 5.5.  The  moduli spaces may have
several components:
Every component  contains a Zariski open subset which is a  holomorphic
$\C^*$-bundle.  For some components,  this $\C^*$-bundle consists only  of
pairs $({\cal E},\varphi)$ with
${\cal E}$ stable as a  bundle; components of this type can be obtained by
compactifying  the
corresponding $\C^*$-bundle with a  Donaldson moduli space  at infinity. In
the other direction,
the component is not compact, but has a {\sl natural compactification}
obtained by adding spaces
associated with   Seiberg-Witten moduli spaces. The  other components can
also be
naturally compactified  by  using   Seiberg-Witten moduli spaces in
both directions.

This compactification process, as well as the corresponding differential
geometric
interpretation will be the subject of a later paper. \footnote{ After
having completed our
results we received a manuscript by Labastida and Marino [LM] in which
related ideas are
proposed from a physical point of view, and physical implications are
discussed }
\section{$Spin^h$-structures}

The quaternionic spinor group is defined as
$$Spin^h:=Spin\times_{\Z/2}Sp(1)=Spin\times_{\Z/2}SU(2)\ ,$$
and   fits in the  exact
sequences
$$\begin{array}{c}
1\map Sp(1)\map Spin^h\stackrel{\pi}{\map}SO\map 1\\
1\map Spin\map Spin^h\stackrel{\delta}{\map}\ PU(2)\map 1\end{array}
\eqno{(1)}$$
These can be combined in the sequence
$$1\map\Z/2\map Spin^h\textmap{(\pi,\delta)} SO\times PU(2)\map 1\eqno{(2)}$$

In dimension 4, $Spin^h(4)$ has a simple description, coming from the splitting
$Spin(4)=SU(2)\times SU(2)$:
$$Spin^h(4)=\qmod{SU(2)\times SU(2) \times SU(2)}{\Z/2}$$
with $\Z/2=\langle(-\id,-\id,-\id)\rangle$. There is another useful way to
think of $Spin^h(4)$:
let $G$ be the group
$$G:=\{(a,b,c)\in U(2)\times U(2)\times U(2)|\ \det a=\det b= \det c\}\ .$$
One has an obvious isomorphism
$Spin^h(4)=\qmod{G}{S^1}\ $ , and a commutative diagram with exact rows
$$\matrix{1\rightarrow&\Z_2&\map &SU(2)\times SU(2)\times SU(2)&\map
&Spin^h(4)&\rightarrow 1\cr
           &\downarrow&  &\downarrow&&\parallel&\cr
          1\rightarrow &S^1&\map &G&\map &Spin^h(4)&\rightarrow 1\cr
          }\eqno{(3)}$$

\begin{dt} Let $P$ be a principal $SO$-bundle over a space $X$. A
$Spin^h$-structure in $P$ is
a pair consisting of a $Spin^h$ bundle $P^h$ and an isomorphism $P\simeq
P^h\times_\pi SO$.
The $PU(2)$-bundle associated with a $Spin^h$-structure  is the bundle
$P^h\times_\delta
PU(2)$.
\end{dt}
\begin{lm} A principal $SO$-bundle admits a $Spin^h$-structure iff there exists
a
$PU(2)$-bundle with the same second Stiefel-Whitney class.
\end{lm}
\pf This follows from the cohomology sequence
$$\map H^1(X,\underline{Spin^h})\map
H^1(X,\underline{{\phantom(}SO{\phantom)}}\times\underline{PU(2)})\textmap{\
beta}
H^2(X,\Z/2)$$
associated  to (2), since the connecting  homomorphism $\beta$ is given by
taking the sum of
the second Stiefel-Whitney classes of the two factors.
\qed

In this paper we will  only use $Spin^h$-structures in $SO(4)$-bundles
whose second Stiefel
Whitney class admit   {\sl integral} lifts. Then we have:
\begin{lm}

Let $P$ be  a principal $SO(4)$-bundle whose second Stiefel-Whitney class
$w_2(P)$ is the
reduction of an integral class.

Isomorphism classes of $Spin^h(4)$-structures in $P$ are in 1-1
correspondence with
equivalence classes of triples consisting of  a $Spin^c(4)$-structure
$\qmod{P^c}{S^1}\simeq P$ in
$P$, a
$U(2)$-bundle $E$, and an isomorphism $\det P^c\simeq\det E$, where two
triples are
equivalent if they can be obtained from each other by tensoring with an
$S^1$-bundle.
\end{lm}
\pf The cohomology sequence associated with the second row in (3) shows that
$Spin^h$-structures in bundles whose second Stiefel-Whitney classes admit
integral lifts are
given by $G$-structures  modulo tensoring with $S^1$-bundles. On the other
hand, to give a
$G$-structure in $P$ simply means to give a triple $(\Sigma^+,\Sigma^-,E)$
of $U(2)$-bundles
together with isomorphisms
$$\det\Sigma^+\simeq\det\Sigma^-\simeq\det E \ .$$

This is equivalent to giving  a triple consisiting of a
$Spin^c(4)$-structure $\qmod{P^c}{S^1}\simeq
P$ in $P$, a $U(2)$-bundle, and an isomorphism $\det P^c\simeq\det E$.
\qed

In the situation of this lemma, we get well defined vector bundles
$${\cal H}^{\pm}:=\Sigma^{\pm}\otimes E^{\vee}\  $$
 depending only on the $Spin^h$-structure  and
not on the chosen $G$-lifting. These spinor bundles have the following
intrinsic interpretation:
identify
$SU(2)\times_{\Z/2} SU(2)$ with $SO(4)$, and denote by
$$\pi_{ij}:Spin^h\map SO(4)$$
 the projections of $Spin^h(4)=\qmod{SU(2)\times SU(2)\times SU(2)}{\Z/2}$
onto the indicated
factors  ($\pi=\pi_{12}$). Using the inclusion $SO(4)\subset SU(4)$, we can
form three
$SU(4)$-vector bundles $P^h\times_{\pi_{ij}}\C^4$,
$(i,j)\in\{(1,2),(1,3),(2,3)\}$.

Under the conditions of the previous lemma  we have
$${\cal H}^+=P^h\times_{\pi_{13}}\C^4\ ,\ \ {\cal
H}^-=P^h\times_{\pi_{23}}\C^4\ ,\ \
\Sigma^+\otimes(\Sigma^-)^{\vee}=P^h\times_{\pi}\C^4  \  .$$

The $PU(2)$-bundle $P^h\times_\delta PU(2)$ associated with the
$Spin^h$-structure
$\qmod{P^c}{S^1}\simeq P$ has in this case a very simple description: it is
the projectivization
$P(E)$ of the
$U(2)$-bundle $E$.

\section{ The quaternionic monopole equations}

Let $(X,g)$ be an oriented Riemannian 4-manifold with orthonormal frame
bundle $P$. The exact
sequence (2) in the previous section shows two things: first, isomorphism
classes of
$PU(2)$-bundles with second Stiefel-Whitney class   equal to $w_2(P)$ are in
1-1 correspondence with orbits of $Spin^h(4)$-structures in $P$ under the
action of
$H^1(X,\Z /2)$; second, $Spin^h(4)$-connections in a $Spin^h(4)$-bundle
$P^h$ which induce the
Levi-Civita connection in  $P$ correspond bijectively to connections in the
associated
$PU(2)$-bundle $P^h\times_\delta PU(2)$.

Now it is well known that $w_2(P)=w_2(X)$ is always the reduction of an
integral class [HH], so
that we can think of a $Spin^h$-structure in $P$ as a triple
$(\Sigma^+,\Sigma^-,E)$ of
$U(2)$-bundles with isomorphisms $\det\Sigma^+\simeq\det\Sigma^-\simeq\det
E$ modulo
tensoring with unitary line bundles. We denote the $Spin^h(4)$-connection
corresponding to a
connection $A\in{\cal A}(P(E))$ in the associated $PU(2)$-bundle  by  $\hat A$.
\begin{re} Given a fixed $U(1)$-connection $c$ in $\det E$, the elements in
${\cal A}(P(E))$ can
be identified with those $U(2)$-connections in $E$, which induce the fixed
connection $c$.
\end{re}

Now view a $Spin^h(4)$-structure in $P$ as a $Spin^c(4)$-structure
$\qmod{P^c}{S^1}\simeq P$
together with a $U(2)$-bundle $E$ and an isomorphism $\det P^c\simeq\det E$
. Recall
that the choice of $\qmod{P^c}{S^1}\simeq P$ induces an isomorphism
$$\gamma:\Lambda^1\otimes\C\map(\Sigma^+)^{\vee}\otimes\Sigma^-$$
which extends to a homomorphism
$$\Lambda^1\otimes\C\map\End_0(\Sigma^+\oplus\Sigma^-)\ ,$$
mapping the bundle $\Lambda^1$ of real 1-forms into the bundle of
trace-free skew-Hermitian
endomorphisms. The induced homomorphism
$$\Gamma:\Lambda^2\otimes\C\map \End_0(\Sigma^+\oplus\Sigma^-)$$
maps the subbundles $\Lambda^2_{\pm}\otimes\C$ isomorphically onto the bundles
$\End_0(\Sigma^{\pm})$, and identifies $\Lambda_{\pm}$ with the trace-free,
skew-Hermitian
endomorphisms ([H], [OT1]).
\begin{dt}
Let $P^h\times_{\pi}SO(4)\simeq P$ be a $\Spin^h(4)$-structure in $P$ with
spinor bundle ${\cal
H}:={\cal H}^+\oplus{\cal H}^-$ and associated $PU(2)$-bundle $P(E)$.
Choose a connection
$A\in{\cal A}(P(E))$, and let $\hat A$ be the corresponding
$Spin^h(4)$-connection  in $P^h$.
The associated Dirac opearor  is defined as the composition
$$\Dr_A:A^0({\cal H})\textmap{\nabla_{\hat A}} A^1({\cal
H})\textmap{\gamma} A^0({\cal
H})\ ,$$
where $\nabla_{\hat A}$ is the covariant derivative of $\hat A$ and
$\gamma$ the Clifford
multiplication.
\end{dt}
Note that the restricted operators
$$\Dr_A : A^0({\cal H}^{\pm})\map A^0({\cal H}^{\mp})$$
interchange the positive and negative half-spinors.

Let $s$ be the scalar curvature of $(X,g)$.
\begin{pr} The Dirac operator $\Dr_A:A^0({\cal H})\map A^0({\cal H})$ is an
elliptic, selfadjoint
operator whose Laplacian satisfies the Weitzenb\"ock formula
$$\Dr_A^2=\nabla_{\hat A}^*\nabla_{\hat A}+\Gamma(F_A)+\frac{s}{4}id_{\cal
H}\eqno{(4)}$$
\end{pr}
\pf Choose a $Spin^c(4)$-structure  $\qmod{P^c}{S^1}\simeq P$ and a
$S^1$-connection $c$ in the
unitary line bundle $\det P^c$. The connection $A\in {\cal A}(P(E))$lifts
to a unique
$U(2)$-connection $C$ in the bundle $E^{\vee}$ which induces the dual
connection of $c$ in $\det
E^{\vee}=(\det P^c)^{\vee}$. In [OT1] we introduced the Dirac operator
$$\Dr_{C,c}:A^0(\Sigma \otimes  E^{\vee}) \map A^0(\Sigma \otimes E^{\vee})\
;$$
by construction it coincides with the operator $\Dr_A:A^0({\cal H})\map
A^0({\cal H})$, and its
Weitzenb\"ock formula reads
$$\Dr_{C,c}^2=\nabla_{\hat A}^*\nabla_{\hat
A}+\Gamma(F_{C,c})+\frac{s}{4}id_{\cal H}\ ,$$
where $F_{C,c}=F_C+\frac{1}{2}F_c\id_{E^{\vee}}\in A^2( \End E^{\vee})$.
Substituting $$F_C
=\frac{1}{2}\tr(F_C)\id_{E^{\vee}}+F_A$$ and using
$\frac{1}{2}\tr(F_C)=-\frac{1}{2}F_c$ we
get the Weitzenb\"ock formula (4) for $ \Dr_A$.

\qed

Consider now a section $ \Psi\in A^0({\cal H}^{\pm})$. We denote by
$$(\Psi\bar\Psi)_0\in A^0(\End_0 \Sigma^{\pm}\otimes\End_0 E^{\vee})$$
the projection of $\Psi\otimes\bar\Psi\in A^0(\End{\cal H}^{\pm})$ onto the
fourth summand in
the decomposition
$$\End({\cal H}^{\pm})=\C\id\oplus\End_0\Sigma^{\pm}\otimes\End_0
E^{\vee}\otimes(\End_0\Sigma^{\pm}\otimes\End_0 E^{\vee})\ .$$
$(\Psi\bar\Psi)_0$ is a Hermitian endomorphism which is trace-free in both
factors.
\begin{dt}Choose a $Spin^h(4)$-structure in $P$ with spinor bundle ${\cal
H}$ and associated
$PU(2)$-bundle $P(E)$. The quaternionic monopole equations for the pair
$(A,\Psi)\in{\cal
A}(P(E))\times A^0({\cal H})$ are the following equations:
$$\left\{\begin{array}{ccc}\Dr_{A}\Psi&=&0\ \ \ \ \\
\Gamma(F_{A}^+)&=&(\Psi\bar\Psi)_0 \ .\ \
\end{array}\right.\eqno{(SW^h) }$$
\end{dt}
The following result  is the analog of Witten's formula  in the
quaternionic case (see [W], \S 3 ):
\begin{pr}   Let $\Psi\in A^0({\cal H}^+)$ be a positive half-spinor, $A\in
{\cal A}(P(E))$ a connection in $P(E)$. Then we have
$$\parallel\Dr_{A}\Psi\parallel^2+
\frac{1}{2}\parallel\Gamma(F_{A}^+)-(\Psi\bar\Psi)_0\parallel^2=$$ $$=
\parallel\nabla_{\hat A}\Psi\parallel^2+
\frac{1}{2}\parallel
F_{A}^+\parallel^2+\frac{1}{2}\parallel(\Psi\bar\Psi)_0\parallel^2+
\frac{1}{4}\int\limits_X s|\Psi|^2. \eqno{(5)}$$
\end{pr}

\pf  The pointwise inner product $(\Gamma(F_A)\Psi,\Psi)$ for a positive
half-spinor
$\Psi$ simplifies:
$(\Gamma(F_A)\Psi,\Psi)= (\Gamma(F_A^+)\Psi,\Psi)=
(\Gamma(F_A^+),(\Psi\bar\Psi)_0)$, since $\Gamma(F_A^-)$ vanishes on
$A^0({\cal H}^+)$,
and since $ \Gamma(F_A^+)$ is trace-free in both arguments.

Using the Weitzenb\"ock formula (5), we find
$$(\Dr_A\Psi,\Psi)=(\nabla_{\hat A}^*\nabla_{\hat
A}\Psi,\Psi)+(\Gamma(F_A^+),(\Psi\bar\Psi)_0)+\frac{s}{4}|\Psi|^2\
,\eqno{(6)}$$
which shows that
$$(\Dr_A^2\Psi,\Psi)+\frac{1}{2}|\Gamma(F_{A}^+)-(\Psi\bar\Psi)_0|^2=
(\nabla_{\hat A}^*\nabla_{\hat A}\Psi,\Psi)+
\frac{1}{2}|F_{A}^+|^2+\frac{1}{2}|(\Psi\bar\Psi)_0|^2+\frac{s}{4}|\Psi|^2
$$
The identity (5) follows  by integration over $X$.

\section{Moduli spaces of quaternionic monopoles}

Let $E$ be $U(2)$-bundle with $w_2(P)\equiv \overline{c_1(E)} $ (mod 2),
and let $c$ be a
fixed $S^1$-connection in $\det E^{\vee}$. We identify ${\cal A}(P(E))$
with the space
${\cal A}_c(E^{\vee})$ of $U(2)$-connections in $E^{\vee}$ which induce the
fixed
connection in $\det E^{\vee}$, and we set:
$${\cal A}:={\cal A}_c(E^{\vee})\times A^0({\cal H}^+)$$
The natural gauge group is the group ${\cal G}$ consisting of unitary
automorphisms in
$E^{\vee}$ which induce the identity in $\det E^{\vee}$. ${\cal G}$ acts on
${\cal A}$ from
the right  in a natural way. Let ${\cal A}^*\subset {\cal A}$ be the open
subset of ${\cal
A}$ consisting of pairs $(C,\Psi)$ whose stabilizer ${\cal G}_{(C,\Psi)}$
is contained in the
center $\Z/2=\{\pm\id_E\}$  of the gauge group.
\begin{re} A pair $(C,\Psi)$ does not belong to ${\cal A}^*$ iff $\Psi=0$
and $C$ is a
reducible connection.
\end{re}
Indeed, the isotropy group of ${\cal G}$ acting only on the first factor
${\cal A}_c(E^{\vee})$ is
the centralizer of the holonomy of $C$ in $SU(2)$. The latter is $S^1$ if
$C$ is reducible,
and $\Z/2$ in the irreducible case.
\qed

A pair belonging to ${\cal A}^*$ will be called irreducible. Note that the
stabilizer of
\underbar{any} pair with vanishing second componenent
$\Psi$ contains
$\Z/2$.

{}From now on we also assume that ${\cal A}$ and ${\cal G}$ are completed
with respect to
suitable Sobolev norms $L^2_k$, such that ${\cal G}$ becomes a Hilbert Lie
group acting
smoothly on ${\cal A}$. Let ${\cal B}:=\qmod{{\cal A}}{{\cal G}}$ be the
quotient, ${\cal
B}^*:=\qmod{{\cal A}^*}{{\cal G}}$, and denote the orbit-map  $[\ ]:{\cal
A}\map{\cal B}$
 by $\pi$.

An element in ${\cal A}^{*}$ will be called {\sl strongly irreducible} if
its stabilizer is
trivial. Let ${\cal A}^{**}\subset{\cal A}^*$ be the subset of   strongly
irreducible  pairs,
and put ${\cal B}^{**}:=\qmod{{\cal A}^{**}}{{\cal G}}$.
\begin{pr} ${\cal B}$ is a Hausdorff space. ${\cal B}^{**}\subset{\cal B}$
is open and has
the structure of a differentiable Hilbert manifold. The map ${\cal
A}^{**}\map{\cal
B}^{**}$ is a differentiable principal ${\cal G}$-bundle.
\end{pr}
\pf Standard, cf. [DK], [FU]. \\

Fix a point $p=(C,\Psi)\in{\cal A}$. The  differential of the map ${\cal
G}\map{\cal A}$
given by the action of ${\cal G}$ on $p$ is the map
$$
\begin{array}{cccc}D^0_{p}&:A^0(su(E^{\vee}))&\map& A^1(su(E^{\vee}))\oplus
A^0(\Sigma^+\otimes E^{\vee})\\ &f&\longmapsto&(D_{C}(f),-f\Psi)
\end{array}$$

Setting
$$N_{p}(\varepsilon):=\{\beta\in A^1(su(E^{\vee}))\oplus
A^0(\Sigma^+\otimes E^{\vee})|\ {D^0_{p}}^*\beta=0,\
||\beta||<\varepsilon\}\ ,$$
for $\varepsilon>0$ sufficiently small, one obtains local slices for the
action of ${\cal G}$
on ${\cal A}^{**}$  and charts
$\pi|_{N_p(\varepsilon)}:N_p(\varepsilon)\map {\cal B}^{**}$ for
${\cal B}^{**}$.
\qed

Note that the curvature $F_A$ of a connection in $P(E)$ equals the
trace-free part $F_C^0$ of the
curvature of the corresponding connection $C\in{\cal A}_c(E^{\vee})$.

Using the identification ${\cal A}(P(E))={\cal A}_c(E^{\vee})$, we can
rewrite the
quaternionic monopole equations in terms of pairs $(C,\Psi)\in{\cal A}$.
Let ${\cal
A}^{SW^h}\subset {\cal A}$ be the space of solutions.

\begin{dt} Fix a $Spin^h$-structure in $P$. The moduli space of quaternionic
monopoles is
the quotient  ${\cal M}:=\qmod{{\cal A}^{SW^h}}{{\cal G}}$.   We denote by
${\cal
M}^{**}:=\qmod{({\cal A}^{SW^h}\cap{\cal A}^{**})}{{\cal G}}$,  ${\cal
M}^{*}:=\qmod{({\cal A}^{SW^h}\cap{\cal A}^{*})}{{\cal G}}$   the subspaces of
(strongly) irreducible monopoles.
\end{dt}

The tangent space to ${\cal A}^{SW^h}$ at $p=(C,\Psi)\in {\cal A}$ is the
kernel of the
operator
$$
\begin{array}{c}D^1_{p} :A^1(su(E^{\vee}))\oplus A^0(\Sigma^+\otimes
E^{\vee}) \map
A^0(su(\Sigma^+)\otimes su(E^{\vee}))\oplus A^0(\Sigma^-\otimes
E^{\vee})\end{array}$$
defined by
$$D^1_{p}((\alpha,\psi))= \left(\Gamma(D_C^+(\alpha))-[(\psi\bar\Psi)_0+
(\Psi\bar\psi)_0],\Dr_{C,c}\psi +\gamma(\alpha)\Psi\right)    \ ,
 $$
where we consider $\gamma(\alpha)$ as map  $\gamma(\alpha):\Sigma^+\map
\Sigma^-\otimes su(E^{\vee})$. Clearly $D^1_p\circ D_p^0=0$, since the
monopole equations
are gauge invariant.

Using the isomorphism $\Gamma^{-1}: A^0(su(\Sigma^+))\map A^2_+$, we can
consider
$D^1_p$ as an operator $D^1_p:A^1(su(E^{\vee}))\oplus A^0(\Sigma^+\otimes
E^{\vee}) \map
A^2_+(su(E^{\vee}))\oplus A^0(\Sigma^-\otimes E^{\vee})$.

Let $\sigma(X)$ and $e(X)$ be the signature and the topological Euler
characteristic of the
oriented manifold $X$.

\begin{pr} For a solution $p=(C,\Psi)\in{\cal A}^{SW^h}$, the complex
$$0\rightarrow A^0 su(E^{\vee}) \textmap{{D}^0_p}
A^1 su(E^{\vee}) \oplus A^0 {\cal H}^+ \textmap{{D}^1_p}A^2_+ su(E^{\vee})
\oplus A^0 {\cal
H}^- \rightarrow 0\eqno{({\cal C}_p)}$$
is elliptic  and its index is
$$\frac{3}{2}(4c_2(E^{\vee})-c_1(E^{\vee})^2)-\frac{1}{2}(3e(X)+4\sigma(X))
\ .\eqno{(7)}$$
\end{pr}
\pf
The complex ${\cal C}_p$ has  the same symbol sequence as
$$0\rightarrow A^0 su(E^{\vee}) \stackrel{(D_C,0)}{\map}
A^1 su(E^{\vee}) \oplus A^0 {\cal
H}^+  \stackrel{(D_C^+,\Dr_{C,c})}{\map}A^2_+ su(E^{\vee}) \oplus A^0 {\cal
H}^- \rightarrow 0 $$
which is an elliptic complex with index
$$2(4c_2(E^{\vee})-c_1(E^{\vee})^2)-\frac{3}{2}(\sigma(X)+e(X))+index
\Dr_{C,c}\  .$$

 The latter term is
$$index\Dr_{C,c}=[ch(E^{\vee})e^{\frac{1}{2} c_1(E^{\vee})}\hat
A(X)]_4=-2c_2(E^{\vee})+\frac{1}{2}c_1(E^{\vee})^2-\frac{1}{2}\sigma(X)\ .$$
\qed
\begin{re} The integer in (7) is always an even number if $X$ admits almost
complex
structures.
\end{re}

Our next step is to endow the spaces ${\cal M}^{**}$ (${\cal M}^{*}$) with
the structure  of
a real analytic  space (orbifold).

In the first case (compare with [FU], [DK], [OT1], [LT]), we have an
analytic map
$\sigma:{\cal A}\map A^2_+(su(E^{\vee}))\oplus A^0( {\cal H}^-)$ defined by
$$\sigma(C,\Psi)=
\left((F_C^0)^+-\Gamma^{-1}(\Psi\bar\Psi)_0,\Dr_{C,c}\Psi\right)$$
 which gives rise to a section $\tilde\sigma$ in the bundle ${\cal
A}^{**}\times_{{\cal G}}
\left(A^2_+(su(E^{\vee}))\oplus A^0( {\cal H}^-) \right)$. We endow ${\cal
M}^{**}$ with a real
analytic structure by identifying it with the vanishing locus
$Z(\tilde\sigma)$ of
$\tilde\sigma$, regarded as a  subspace of the Hilbert manifold ${\cal
B}^{**}$  (in Douady's sense) ([M], [LT]).

Now fix a point $p=(C,\Psi)\in{\cal A}^*$.  We put
$$S_p(\varepsilon):=\{p+\beta|\ \beta\in A^1su(E^{\vee}) \oplus A^0 {\cal
H}^+ ,\
D^0_p{D^0_p}^* \beta +{D^1}^*_p\sigma(p+\beta)=0,\ ||\beta||<\varepsilon\}\ .$$
\begin{cl} For sufficiently small $\varepsilon>0$, $S_p(\varepsilon)$  is a
finite
dimensional submanifold of ${\cal A}$ which is contained in the slice
$N_p(\varepsilon)$
and whose tangent space at $p$ is the first harmonic space $\H^1_p$ of the
deformation
complex ${\cal C}_p$.
\end{cl}

To prove this claim, we consider the map
$$s_p:A^1(su(E^{\vee}))\oplus A^0({\cal H}^+)\map \im( D^0_p)\oplus\im
(D^1_p)^*$$
given by the left hand terms in the equations defining $S_p(\varepsilon)$.
The derivative of
$s_p$ at 0 is the first Laplacian
$$\Delta^1_p:A^1(su(E^{\vee}))\oplus A^0({\cal H}^+)\map \im(
D^0_p)\oplus\im (D^1_p)^*$$
associated with the elliptic complex ${\cal C}_p$, hence $s_p$ is  a
submersion in 0. This
proves the claim.
\qed

The intersection ${\cal A}^{SW^h}\cap N_p(\varepsilon)=Z(\sigma)\cap
N_p(\varepsilon)$
of the space of solutions with the standard slice through $p$ is contained in
$S_p(\varepsilon)$ and can be identified with the finite dimensional model
$$Z(\sigma)\cap N_p(\varepsilon)=Z(\sigma|_{S_p(\varepsilon)})\ .$$
If $p\in {\cal A}^{**}$ is strongly irreducible, then the map
$$\pi|_{Z(\sigma|_{S_p(\varepsilon)})}:Z(\sigma|_{S_p(\varepsilon)})\map
{\cal M}^{**}$$
is a local parametrization of ${\cal M}^*$ at $p$, hence
$Z(\sigma|_{S_p(\varepsilon)})$ is
a local model for the moduli space around $p$.

If $p\in{\cal A}^*\setminus{\cal A}^{**}$ is irreducible but not strongly
irreducible, then
necessarily $\Psi=0$, and the isotropy group  ${\cal G}_p=\Z/2$ acts on
$S_p(\varepsilon)$. Since
$\sigma$ is $\Z/2$-equivariant, we obtain an induced action on
$Z(\sigma|_{S_p(\varepsilon)})$. In
this case $\pi|_{Z(\sigma|_{S_p(\varepsilon)})}$ induces a homeomorphism of
the quotient
$\qmod{Z(\sigma|_{S_p(\varepsilon)})}{\Z/2}$ with an open neighbourhood of
$p$ in ${\cal
M}^*$, and ${\cal M}^*$ becomes an orbifold at $p$, if we use the map
$$\pi|_{Z(\sigma|_{S_p})}:Z(\sigma|_{S_p(\varepsilon)})\map {\cal M}^*$$
as an orbifold chart.

\begin{re} Using a real analytic isomorphism which identifies the germ of
$S_p(\varepsilon)$ at $p$ with the germ of $\H_p^1=T_p(S_p(\varepsilon))$
at 0, we
obtain a local model of Kuranishi-type for ${\cal M}^*$ at $p$.
\end{re}
\begin{re}  The points in ${\cal D}^*:={\cal M}^*\setminus{\cal M}^{**}$
have the form
$[(C,0)]$, where
$C$ is projectively anti-self-dual, i.e $(F_C^0)^+=0$. There is a natural
finite map
$${\cal D}^*\map{\cal M}(P(E^{\vee}))$$
into the Donaldson moduli space of $PU(2)$-instantons in
$P(E^{\vee})$, which maps ${\cal D}^*$ isomorphically  onto ${\cal
M}(P(E^{\vee})^*$ if
$H^1(X,\Z/2)=0$. In general   ${\cal D}^*$ and ${\cal M}(P(E^{\vee})^*$
cannot be identified.
The difference comes from the fact that our gauge group is $SU(E^{\vee})$,
whereas the
$PU(2)$-instantons are classified modulo $PU(E^{\vee})$.
\end{re}

For simplicity we shall however refer to ${\cal D}^*$ as  Donaldson
instanton moduli space.

Concluding, we get

\begin{pr} ${\cal M}^{**}$ is a real analytic space. ${\cal M}^*$ is a real
analytic orbi\-fold,
and the points in  ${\cal M}^*\setminus{\cal M}^{**}$ have neighbourhoods
modeled on
$\Z/2$-quotients. ${\cal M}^*\setminus{\cal M}^{**}$ can be identified as a
set with the
Donaldson moduli space ${\cal D}^*$ of irreducible projectively
anti-self-dual connections in
$E^{\vee}$ with fixed determinant $c$.
\end{pr}

The local structure of the moduli space ${\cal M}$ in   reducible points, which
correspond to  pairs formed by a reducible instanton and a trivial spinor,
can also be
described using the method above  (compare with [DK]).

Let $\MM^{SW}\subset{\cal M}$ be the subspace of ${\cal M}$ consisting of
all  orbits of the
form
$(C,\Psi)\cdot SU(E^{\vee})$, where $C$ is a reducible connection  and
$\Psi$ belongs to one
of the summands. Let $L:=\det\Sigma^{\pm}=\det E$. It is easy to see that
$$\MM^{SW}\simeq\union\limits_{\matrix{^{S {\rm \ summand}}\cr ^{{\rm of} \
E^{\vee}}\cr}}{\cal M}^{SW}_{L \otimes S^{\otimes 2}} \ ,$$
where  ${\cal M}^{SW}_M$ denotes  the rank-1 Seiberg-Witten moduli space
associated to a $Spin^c(4)$-structure of determinant $M$.

 The  fact that the moduli spaces of quaternionic monopoles contain
Donaldson moduli spaces  as well of Seiberg-Witten moduli spaces suggests
that they could
provide a method for comparing the invariants given by the two theories.

\section{Quaternionic monopoles on K\"ahler surfaces}

Let $(X,g)$ be a K\"ahler surface with canonical $Spin^c(4)$-structure; in
this case
$\Sigma^+=\Lambda^{00}\oplus\Lambda^{02}$, and $\Sigma^-=\Lambda^{01}$. A
$Spin^h(4)$-structure  in the frame bundle is given by a unitary vector
bundle $E$
together with an isomorphism $\det E\simeq \Lambda^{02}$. A
$Spin^h(4)$-connection
$\hat A$ corresponds to a $PU(2)$-connection $A$ in the associated bundle
$P(E)$, or
alternatively, to a unitary connection $C$ in $E^{\vee}$ which induces a fixed
$S^1$-connection $c$ in $\Lambda^{20}$.  Recall that the curvature $F_A$ of
$A$ equals the
trace-free component $F_C^0$ of $F_C$.

If we choose
$c$ to be the Chern connection in the canonical  bundle $\Lambda^{20}=K_X$,
then the
$Spin^h(4)$-connection in
${\cal H}=\Sigma \otimes E^{\vee}$ is simply the tensor product of  the
canonical
connection in $ \Sigma=\Sigma^{+}\oplus\Sigma^-$ and the connection $C$.

A positive quaternionic spinor $\Psi\in A^0 ( {\cal H}^+) $ can be written as
$\Psi=\varphi+\alpha$, with $\varphi\in A^0( E^{\vee}) $, and $\alpha\in
A^{02}( E^{\vee})$.

\begin{pr} Let $C$ be a unitary connection in $E^{\vee}$ inducing the Chern
connection $c$ in
$\det E^{\vee}=K_X$. A pair $(C,\varphi+\alpha)$ solves the quaternionic
monopole
equations if and only if $F_C$ is of type $(1,1)$
and one  of the following  conditions holds
$$\matrix{1.\ \alpha=0\ ,\ \bar\partial_C\varphi=0\ and\   i\Lambda_g
F_{C}^0 +\frac{1}{2} (\varphi\otimes\bar\varphi)_0\ =\ 0\ , \cr
  2.  \ \varphi=0\ ,\ \partial_C\alpha=0\ and\   i\Lambda_g
F_{C}^0 -\frac{1}{2}* (\alpha\otimes\bar\alpha)_0=0\ .\cr}\eqno{(8)}$$
\end{pr}

\pf Using the notation  in the proof of the Weitzenb\"ock formula, we have
$F_{C,c}=\frac{1}{2}(\tr F_C+F_c)\id_{E^{\vee}}+ F_A=F_A=F_C^0\in A^2(su
(E^{\vee}))$. By
Proposition 2.6 of [OT1]  the quaternionic Seiberg-Witten equations become
$$\left\{
\begin{array}{lll}
   F_{A}^{20}&=&-\frac{1}{2}(\varphi\otimes\bar\alpha)_0\\
  F_{A}^{02}&=&\frac{1}{2}(\alpha\otimes\bar\varphi)_0\\
  i\Lambda_g F_{A}&=&-\frac{1}{2}\left[(\varphi\otimes\bar\varphi)_0-
  *(\alpha\otimes\bar\alpha)_0\right]\\
  \bar\partial_C\varphi&=&i\Lambda_g\partial_C\alpha\ \ \ \ \ \ \ \ \ \ \ \
\ \ \ \ \ \ \
.\end{array}\right.$$

 Note that the right-hand side of formula (5) is invariant under Witten's
transformation
$(C,\varphi+\alpha)\longmapsto (C, \varphi-\alpha)$. Therefore, every
solution satisfies
$F_A^{20}= F_A^{02}=0$,
and $(\varphi\otimes\bar\alpha)_0=(\alpha\otimes\bar\varphi)_0=0$. Elementary
computations show that this can happen only if $\varphi=0$ or $\alpha=0$.
On the other
hand, since the Chern connection in $K_X$ is integrable, we also get
$F_C^{20}=F_C^{02}=0$.
\qed
\begin{re} The second case in this proposition reduces to the first: in fact,
if
$\varphi=0$ and $\alpha\in A^{02}(E^{\vee})$ satisfies
$i\Lambda_g\partial\alpha=0$, we set
$\psi:=\bar\alpha\in A^{20}(\bar E^{\vee})=A^0(\Lambda^{20}\otimes
E)=A^0(E^{\vee})$, and
we obtain
$\bar\partial_C\psi=\overline{\partial_C\bar\psi}=\overline{\partial_C\alpha
}=0$. Here
we used the fact that $\Lambda_g:\Lambda^{12}\map \Lambda^{01}$ is an
isomorphism, the
adjoint of the Lefschetz isomorphism $\cdot\wedge\omega_g$ [LT]. A simple
calculation
in coordinates gives
$-*(\alpha\otimes\bar\alpha)_0=(\bar\alpha\otimes\overline{\bar\alpha})_0=
(\psi\otimes\bar\psi)_0$.
\end{re}

\section{Stability}

Let $(X,g)$ be a compact K\"ahler manifold of arbitrary dimension, $E$ a
differentiable
vector bundle, and let ${\cal L}$ be a fixed holomorphic line
bundle, whose underlying differentiable line bundle is $L:=\det E$.
\begin{dt} A holomorphic pair of type $(E,{\cal L})$ is a pair $({\cal
E},\varphi)$
consisting of a holomorphic bundle ${\cal E}$ and a section $\varphi\in
H^0(X,{\cal E})$
such that the underlying differentiable bundle of ${\cal E}$ is $E$ and
$\det{\cal E}={\cal
L}$.
\end{dt}
Note that the determinant of the holomorphic bundle ${\cal E}$ is fixed,
not only its
isomorphism type.

Two pairs $({\cal E}_i,\varphi_i)$, $i=1,  2$ of the same type are
isomorphic if there
exists an isomorphism  $f:{\cal E}_1\map{\cal E}_2$ with
$f^*(\varphi_2)=\varphi_1$ and $\det f=\id_{\cal L}$.

In other words, $({\cal E}_i,\varphi_i)$ are isomorphic iff there exists a
complex gauge
transformation $f\in SL(E)$ with $f^*(\varphi_2)=\varphi_1$ such that  $f$ is
holomorphic as a map $f:{\cal E}_1\map{\cal E}_2$.
\begin{dt}
 A holomorphic pair $({\cal E},\varphi)$ is simple if any automorphism  of
it  is
of the form $f=\varepsilon\id_{\cal E}$, where $\varepsilon^ {\rk{\cal E}}=1$.
A pair  $({\cal E},\varphi)$ is strongly simple if its only automorphism
is  $\id_{\cal
E}$.
\end{dt}

Note that a simple pair $({\cal E},\varphi)$ with $\varphi\ne 0$ is stongly
simple,
whereas a pair $({\cal E},0)$ is simple iff ${\cal E}$ is a simple bundle.

Note also that $({\cal E},\varphi)$ is simple iff any trace-free
holomorphic endomorphism
$f$ of ${\cal E}$ with $f(\varphi)=0$ vanishes.

For a nontrivial torsion free sheaf ${\cal F}$ on $X$, we denote by
$\mu_g({\cal F})$ its slope
with respect to the K\"ahler metric $g$. Given a holomorphic bundle ${\cal
E}$ over $X$ and
a holomorphic section $\varphi\in H^0(X,{\cal E})$, we let ${\cal S}({\cal
E})$ be the set of
reflexive subsheaves ${\cal F}\subset{\cal E}$ with $0<\rk({\cal
F})<\rk({\cal E})$, and
we define
$${\cal S}_\varphi({\cal E}):=\{{\cal F}\in{\cal S}({\cal E})|\ \varphi\in
H^0(X,{\cal F})\} \
.$$
Recall the following stability concepts [B2]:
\begin{dt}\hfill{\break}
1. ${\cal E}$ is $\varphi$-stable if
$$\max\left(\mu_g({\cal E}),\sup\limits_{{\cal F}'\in{\cal S}({\cal E})}
\mu_g({\cal
F}')\right)<
\inf\limits_{{\cal F}\in {\cal S}_\varphi({\cal E})}\mu_g(\qmod{{\cal
E}}{{\cal F} })\ .$$
2. Let $\lambda\in\R$ be  a real parameter. The pair $({\cal E},\varphi)$ is
$\lambda$-stable iff
$$\max\left(\mu_g({\cal E}),\sup\limits_{{\cal F}'\in{\cal S}({\cal E})}
\mu_g({\cal
F}')\right)<\lambda<
\inf\limits_{{\cal F}\in {\cal S}_\varphi({\cal E})}\mu_g(\qmod{{\cal
E}}{{\cal F}})\ .$$
3.  $({\cal E},\varphi)$ is called $\lambda$-polystable if ${\cal E}$
splits holomorphically
as ${\cal E}={\cal E}'\oplus{\cal E}''$, such that $\varphi\in H^0(X,{\cal
E}')$, $({\cal
E}',\varphi)$ is a $\lambda$-stable pair, and ${\cal E}''$ is a polystable
vector bundle of
slope $\lambda$.
\end{dt}
{}From now on we restrict ourselves to the case $\rk({\cal E})=2$.
\begin{dt} \hfill{\break}
1.  A holomorphic pair $({\cal E},\varphi)$ of type $({ E},{\cal L})$
is called stable if  one of the following  conditions is satisfied:\\
i) ${\cal E}$ is $\varphi$-stable.\\
ii) $\varphi\ne 0$ and ${\cal E}$ splits in direct sum of line bundle
${\cal E}={\cal
E}'\oplus{\cal E}''$, such that $\varphi\in H^0({\cal E}')$  and the pair
$({\cal
E}',\varphi)$ is $\mu_g({ E})$-stable.\\
2.  A holomorphic pair $({\cal E},\varphi)$ of type $({ E},{\cal L})$
is   called polystable if  it is stable, or $\varphi=0$  and ${\cal E}$ is
a polystable
bundle.
\end{dt}

Note that there is  \underbar{no} parameter $\lambda$ in the stability
concept for
holomorphic pairs of a fixed type. The conditions depend only on the metric
$g$ and on the
slope $\mu_g(E)$ of the underlying differentiable bundle $E$.
\begin{lm} Let $({\cal E},\varphi)$ be a holomorphic pair of type $(E,{\cal
L})$ with
$\varphi\ne 0$. There exists a uniquely determined effective divisor
$D=D_\varphi$ and a
commutative diagram
$$\begin{array}{cccc}

0\map&{\cal O}_X(D)&\textmap{\hat\varphi}&{\cal E}\map {\cal L}(-D)\otimes
J_Z\map 0\
,\\
 &{\scriptstyle D\cdot}\uparrow{\phantom i}&{\scriptstyle\varphi}\nearrow&\\
 &\ {\cal O}_X&\end{array}\eqno{(9)}$$
 with a local complete intersection $Z\subset X$ of codimension 2. The pair
$({\cal
E},\varphi)$ is stable if and only if  $\mu_g({\cal O}_X(D))<\mu_g(E)$.
\end{lm}
\pf $D=D_\varphi$ is the divisorial component of the zero locus
$Z(\varphi)$ of ${\cal E}$
which is defined by the ideal $\im(\varphi^{\vee}:{\cal E}^{\vee}\map{\cal
O}_X)$, and
$\hat\varphi$ is the induced map. The set ${\cal S}_\varphi({\cal E})$
consists precisely
of the line bundles ${\cal F}\subset{\cal O}_X({D})$, so that
$$\inf\limits_{{\cal F}\in
{\cal S}_\varphi({\cal E})}\mu_g(\qmod{{\cal E}}{{\cal
F}})=2\mu_g(E)-\mu_g({\cal
O}_X(D))\ .$$

Suppose $({\cal E},\varphi)$ is stable. If ${\cal E}$ is $\varphi$-stable,
we have
$\mu_g(E)<2\mu_g(E)-\mu_g({\cal O}_X(D))$, which gives the required
inequality. If ${\cal
E}$ is not $\varphi$-stable, then $Z=\emptyset$, the extension (9) splits,
and the
pair $({\cal O}_X(D),\varphi)$ is $\mu_g(E)$-stable, i.e. $\mu_g({\cal
O}_X(D))<\mu_g(E)$.

Conversely, suppose $\mu_g({\cal O}_X(D))<\mu_g(E)$, and assume first that
the extension
(9) does not split. In this case ${\cal E}$ is $\varphi$-stable: in fact,
if ${\cal
F}'\subset {\cal E}$ is an arbitrary line bundle, either ${\cal
F}'\subset{\cal O}_X(D)$, or
the induced map ${\cal F}'\subset{\cal E}\map {\cal J}_Z\otimes{\cal L}(-D)$ is
non-trivial. But then ${\cal F}'\simeq {\cal L}\otimes{\cal
O}_X(-D-\Delta)$ for an
effective divisor $\Delta$ containing $Z$, and we find
$$\mu_g({\cal
F}')=2\mu_g(E)-\mu_g(D)-\mu_g(\Delta)\leq 2\mu_g(E)-\mu_g({\cal
O}_X(D)) \ . $$
 Furthermore, strict inequality holds, unless
$Z=\emptyset$ and the extension (9)  splits, which it  does not by assumption.

In the case of a  split extension, we only have  to notice  that  a pair
$({\cal E}',\varphi)$
 is $\lambda$-stable  for any parameter
$\lambda>\mu_g({\cal E}') $ [B1].
\qed
\begin{re} Consider a pair $({\cal E},\varphi)$ of type $(E,{\cal L})$ with
$\varphi\ne 0$
and associated extension (9). ${\cal E}$ is $\varphi$-stable iff $\mu_g({\cal
O}_X(D))<\mu_g(E)$, and the extension does not split.
\end{re}

Indeed, if the extension splits, then ${\cal E}$ is not $\varphi$-stable, since
$$\mu_g({\cal L}(-D))=\inf\limits_{{\cal F}\in
{\cal S}_\varphi({\cal E})}\mu_g(\qmod{{\cal E}}{{\cal F}}) \ .$$

\section{The projective vortex equation}

Let $E$ be a differentiable vector bundle over a compact K\"ahler manifold
$(X,g)$. We
fix a holomorphic line bundle ${\cal L}$ and a Hermitian metric $l$ in
${\cal L}$. Let
$({\cal E},\varphi)$ be a holomorphic pair of type $(E,{\cal L})$.
\begin{dt} A Hermitian metric in ${\cal E}$ with $\det h=l$ is a solution
of the projective
vortex equation iff the trace free part $F^0_h$ of the curvature $F_h$
satisfies the
equation
$$i\Lambda_g F_h^0 +\frac{1}{2}(\varphi\bar\varphi^h)_0=0\ .\eqno{(V)}$$
\end{dt}
\begin{th} Let $({\cal E},\varphi)$ be a holomorphic pair of type $(E,{\cal
L})$ with
$\rk({\cal E})=2$. Fix a Hermitian metric $l$ in ${\cal L}$.

The pair $({\cal E},\varphi)$ is polystable iff ${\cal E}$ admits a
Hermitian metric $h$
with $\det h=l$ which is a solution of  the projective vortex equation. If
$({\cal
E},\varphi)$ is stable, then the metric $h$ is unique.
\end{th}
\pf
Suppose first that  $h$ is  a solution of the projective vortex equation
$(V)$.  Then we have
$$i\Lambda F_h+\frac{1}{2}(\varphi\bar\varphi^h)=\frac{1}{2}(i\Lambda\tr
F_h+\frac{1}{2}|\varphi|^2)\id_E \ ,
$$
i.e.  $h$ satisfies the weak vortex equation $(V_t)$ associated to the real
function
$t:=\frac{1}{2}(2i\Lambda\tr F_h+|\varphi|^2)$. Therefore, by [OT1],  the
pair $({\cal
E},\varphi)$ is $\lambda$-polystable for the parameter
$\lambda=\frac{(n-1)!}{4\pi}\int\limits_X t vol_g=\mu_g({\cal
E})+\frac{(n-1)!}{8\pi}||\varphi||^2$.

Let $A$ be the Chern connection of $h$, and denote by ${\cal E}'$   the minimal
$A$-invariant
 subbundle which contains $\varphi$. If ${\cal E}'={\cal E}$, then ${\cal E}$
is
$\varphi$-stable  and the  pair $({\cal E},\varphi)$ is stable.

If ${\cal E}'=0$, hence $\varphi=0$, then $h$ is a weak Hermitian-Einstein
metric,
${\cal E}$ is a polystable bundle, and the pair $({\cal E},\varphi)$ is
polystable by
definition.

In the remaining case ${\cal E}'$ is a line bundle and $\varphi\ne 0$.
Let  ${\cal E}'':={\cal E}'^{\bot}$ be the orthogonal complement of ${\cal
E}'$, and let $h'$
and $h''$ be the induced metrics in ${\cal E}'$ and ${\cal E}''$.  We put
$s:=i\Lambda_g\tr
F_h$. Then, since $h=h'\oplus h''$, the projective vortex equation can be
rewritten as:
$$
\left\{\begin{array}{ll}
i\Lambda F_{h'}+\frac{1}{2}(\varphi\bar\varphi^{h'})=&\frac{1}{2}(s+
\frac{1}{2}|\varphi|_{h'}^2)\id_{{\cal E}'}\\
i\Lambda F_{h''}=&\frac{1}{2}(s+\frac{1}{2}|\varphi|_{h'}^2)
\id_{{\cal E}''} \ .
\end{array}\right. $$

The first of these equations   is equivalent to
$$i\Lambda
F_{h'}+\frac{1}{4}(\varphi\bar\varphi^{h'})=\frac{s}{2}\id_{{\cal E}'}\ ,$$
which implies that $({\cal E}',\frac{\varphi}{\sqrt 2})$ is $\mu_g({\cal
E})$-stable by
[OT1].

Conversely, suppose first that   $({\cal E},\varphi)$ is stable. We have to
consider two
cases:\\
\underbar{Case 1}:  ${\cal E}$ is $\varphi$ stable.

Using  Bradlow's existence theorem, we obtain Hermitian metrics in ${\cal
E}$ satisfying
the usual vortex equations associated with suitable chosen $\lambda$ , and,
of course
these metrics all satisfy the equation $(V)$. The problem
is, however, to find a solution   with an a priori given determinant
$l$.

In order to achieve this stronger result, Bradlow's proof has to modified
slightly at some
points:

One starts by fixing  a background metric $k$ such that $\det k=l$. Denote
by $S_0(k)$
the space of   {trace-free} $k$-Hermitian endomorphisms  of $E$, and let
  $\MM et(l)$ be the space of Hermitian metrics in $E$ with $\det h=l$. On
$$ \MM et(l)^p_2:=\{ke^s  |\ s\in L^p_2(S_0(k))\}$$
we define the functional $M_\varphi:\MM et(l)^p_2\map\R$ by
$$M_\varphi(h):=M_D(k,h)+||\varphi||^2_h-||\varphi||^2_k \ .$$
Here $M_D$ is the Donaldson functional, which is known to satisfy   the
identity
$ \frac{d}{dt}M_D(k,h(t))=2\int\limits_X\tr[ h^{-1}(t)\dot h(t) i\Lambda_g
 F_h  ]\ $
for any smooth path of metrics $h(t)$ [Do], [Ko].  Since     $h^{-1}(t)\dot
h(t)$
 is trace-free for a path in $\MM et(l)$, we obtain
$$\frac{d}{dt}|M_D(k,h(t))= 2\int\limits_X\tr[ h^{-1}\dot h(t) i\Lambda_g
 F_h^0]\ .$$
Similarly, for a path of the form $h(t)=he^{ts}$, with $s \in S_0(h)$, we get
$$\frac{d}{dt}_{|_{t=0}}||\varphi||^2_{h_t}=\frac{d}{dt}_{|_{t=0}}\langle
e^{ts}\varphi,\varphi\rangle_h=\left\langle
\frac{d}{dt}_{|_{t=0}}e^{ts}\varphi,\varphi\right\rangle_h=
\left\langle  s ,\varphi\bar\varphi^h\right\rangle_h=\int\limits_X\tr[s
(\varphi\bar\varphi^h)_0] \ .$$
 This means that, putting $m_\varphi(h):=i\Lambda
F_h^0+\frac{1}{2}(\varphi\bar\varphi^h)_0$, we always have
$$\frac{d}{dt}_{|_{t=0}}M_\varphi(he^{ts})= 2\int\limits_X\tr[ s\
m_\varphi(he^{ts})]\ ,\ $$
so that solving the  projective vortex equation  is equivalent  to finding
a critical
point of the functional $M_\varphi$ (compare with Lemma 3.3 [B2]).
\begin{cl}  Suppose $({\cal E},\varphi)$ is simple. Choose $B>0$ and put
$$ \MM et(l)^p_2(B):=\{ h\in \MM et(l)^p_2|\ ||m_\varphi(h)||_{L^p}\leq B\}
\ .$$
Then any
$h\in\MM et(l)^p_2(B)$ which minimizes $M_\varphi$ on $\MM et(l)^p_2(B)$ is
a weak
solution of the projective vortex equation.
\end{cl}
The essential point  is the injectivity
of the operator $s\longmapsto\Delta_h'(s)+\frac{1}{2}[(\varphi\bar\varphi)
s]_0$ acting on
$L^p_2 S_0(h)$. But from
$$\left\langle\Delta_h's+\frac{1}{2}[\varphi\bar\varphi^h) s]_0,\
s\right\rangle_h=||\bar\partial_h(s)||_{h}^2+||s\varphi||^2_{h} \  $$
we see that this operator is injective on trace-free
endomorphisms if $({\cal E},\varphi)$ is simple.
\qed

Now we can follow again   Bradlow's  proof : if  ${\cal E}$ is
$\varphi$-stable, then there exist positive constants $C_1$, $C_2$ such
that for all $s\in
L^p_2 S_0(k) $ with $ke^s\in\MM et(l)^p_2(B)$ the following "main estimate"
holds:
$$\sup|s|\leq  C_1 M_\varphi(ke^{s})+C_2\ .$$
This follows by  applying Proposition 3.2 of [B2] to an arbitrary
$\tau\in\R$ with
$$\max\left(\mu_g({\cal E}),\sup\limits_{\FF'\in\SS(\EE)}\mu_g({\cal
F'})\right)<\frac{(n-1)!\tau
Vol_g(X)}{4\pi}<\inf\limits_{\FF\in\SS_\varphi(\EE)}
\mu_g(\qmod{\cal E}{\cal F})\ ,
$$
since  Bradlow's functional ${\cal
M}_{\varphi,\tau}$  coincides on $\MM et(l)$ with  $M_{\varphi}$ .

It remains to be shown that the existence of this main
estimate implies the existence of a solution of  the projective vortex
equation.

The main estimate implies   that for any $c>0$, the set

$$\{s\in\ L^p_2S_0(k)  |\ ke^s\in \MM et(l)^p_2(B) ,\ \ M_\varphi(ke^{s})<c\}$$
is bounded  in $L^p_2$. Let $(s_i) $ be a sequence in $L^p_2 S_0(k) $ such that
$ke^{s_i}\in \MM et(l)^p_2(B)$ is a minimizing sequence for $M_\varphi$,
and let $s$ be weak
limit. Then $h:=ke^s$ is a weak solution of the projective vortex equation,
which is smooth by
elliptic regularity [B2].

Finally, we have to treat

\underbar{Case 2}: $\varphi\ne 0$, ${\cal E}={\cal E}'\oplus{\cal E}''$,
with $\varphi\in
H^0({\cal E}')$, and the pair $({\cal E}',\varphi)$ is $\mu_g(E)$-stable.\\

We wish  to find
metrics $h'$  and $h''$ in ${\cal E}'$ and ${\cal E}''$, such that  for
$s:=i\Lambda F_l$  the
following equations are satisfied:
$$\left\{ \begin{array}{lll}h'\cdot h''&=&l\\
i\Lambda
F_{h'}+\frac{1}{4}(\varphi\bar\varphi^{h'})&=&\frac{1}{2}s\id_{{\cal E}'}
\\
i\Lambda F_{h''}&=&\frac{1}{2}(s+\frac{1}{2}|\varphi|_{h'}^2)
\id_{{\cal E}''} \ .
\end{array}\right.\ $$
Since the pair $({\cal E}',\frac{1}{\sqrt 2}\varphi)$  is
$\mu_g(E)$-stable, there exists by
[OT1]    a unique Hermitian metric
$h'$ in ${\cal E}'$ solving  the second of these equations. With this
solution   the third
equation   can be rewritten as
$$i\Lambda_g F_{h''}=s-i\Lambda_g F_{h'}\ . $$
Since $\int\limits_X(s-i\Lambda_g F_{h'})=\deg({\cal E}'')$, we can solve
this weak
Hermitian-Einstein equation    by a metric $h''$, which is  unique up to
constant
rescaling. The pro\-duct $h'\cdot h''$ is a metric in
${\cal E}'\otimes{\cal E}''={\cal L}$ which  has the same mean curvature
$s$ as $l$, and therefore
differs from $l$ by a constant factor.  We can now simply rescale $h''$ by
the inverse of this
constant, and we get a pair of metrics satisfying the three equations above.

\qed

\section{Moduli spaces of pairs}

Let $E$ be a differentiable vector bundle of rank $r$ over a K\"ahler
manifold $(X,g)$, and let
${\cal L}$ be a holomorphic line bundle whose underlying differentiable
bundle is $L:=\det E$.
\begin{pr} There exists a  possibly non-Hausdorff  complex analytic
orbifold   ${\cal M}^s(E,{\cal
L} )$ parametrizing isomorphism classes of simple holomorphic pairs of type
$(E,{\cal L})$. The open subset ${\cal M}^{ss}(E,{\cal L})\subset {\cal
M}^{s}(E,{\cal L})$ consisting
of strongly simple pairs is a complex analytic space, and the points in
${\cal M}^s(E,{\cal
L})\setminus{\cal M}^{ss}(E,{\cal L})$ have neighbourhoods modeled on
$\Z/r$-quotients.
\end{pr}

\pf Since we use the same method as in the proof of Proposition 3.9, we
only sketch  the main
ideas.

Let $\bar\lambda$ be the semiconnection defining the holomorphic structure
of ${\cal
L}$, and put $\bar{\cal A}:=\bar{\cal A}_{\bar\lambda}(E)\times A^0(E)$,
where $\bar{\cal
A}_{\bar\lambda}(E)$ denotes the affine space of semiconnections in $E$
inducing $\bar\lambda$
in
$L=\det E$.  The complex gauge group $SL(E)$ acts on $\bar{\cal A}$, and
we write $\bar{\cal A}^s$ ($\bar{\cal A}^{ss})$ for the open subset of
pairs whose stabilizer  is
contained in the center $\Z/r$ of $SL(E)$ ( is trivial).     After suitable
Sobolev
completions,
$\bar{\cal A}^{ss}$ becomes the total space of a holomorphic Hilbert principal
$SL(E)$-bundle over $\bar{\cal B}^{ss}:=\qmod{\bar{\cal A}^{ss}}{SL(E)}$.

A point $(\bar\delta,\varphi)\in \bar{\cal A}$ defines a pair of type
$(E,{\cal L})$ iff  it is
integrable, i.e.    iff it satisfies the following equations:
$$\left\{\begin{array}{lll}F^{02}_{\bar\delta}&=&0\\ \bar\delta\varphi&=&0\ .
 \end{array}\right.\eqno{(10)}$$
Here $F^{02}_{\bar\delta}:=\bar\delta^2$ is a $(0,2)$-form with values in
the bundle $\End_0(E)$
of trace-free endomorphisms.  Moreover, isomorphy of pairs of type
$(E,{\cal L})$ corresponds
to equivalence modulo the action of the complex gauge group $SL(E)$.

Let $\bar\sigma$ be the map ${\cal A}\map A^{02}(\End_0(E))\oplus
A^{01}(E)$ sending a pair
$(\bar\delta,\varphi)$ to the left hand sides of (10).  We endow the sets
${\cal
M}^{ss}_X(E,{\cal L})=\qmod{Z(\sigma)\cap\bar{\cal A}^{ss}}{SL(E)}$   (
${\cal M}^{s }_X(E,{\cal
L})=\qmod{Z(\sigma)\cap\bar{\cal A}^{s }}{SL(E)}$ ) with the structure  of
a complex
analytic space (orbifold) as follows:

${\cal M}^{ss}_X(E,{\cal L})$ is defined to be the vanishing locus of the
section
$\tilde{\bar\sigma}$  in the Hilbert vector bundle
 $\bar{\cal A}^{ss}\times_{SL(E)}\left(A^{02}\End_0(E)\oplus A^{01}E\right)$
over $\bar{\cal B}^{ss}$ which is defined by $\bar\sigma$.

To define the orbifold structure in ${\cal M}^{s }_X(E,{\cal L})$, we use
local models derived from a deformation complex:

Let $\bar p=(\bar\delta,\varphi)\in\bar{\cal A}$ an integrable point. The
associated
 {deformation} {complex} $\bar{\cal D}_{\bar p}$ is the cone over the
evaluation map
$ev^*_\varphi$:
$$ ev^q_\varphi:A^{0q}(\End_0(E))\map A^{0q}(E)\ ,$$
and has the form
$$\matrix{0\rightarrow A^0 (\End_0(E))\textmap{\bar D^0_{\bar p}}
A^{01}(\End_0(E))\oplus
A^0(E)\textmap{\bar D^1_{\bar p}}\cr\  \ \  \  \ \ \textmap{\bar D^1_{\bar p}}
A^{02}(\End_0(E))\oplus A^{01}(E)\textmap{\bar D^2_{\bar
p}}\dots\cr}\eqno{(\bar{\cal D}_{\bar
p})}$$
(compare with [OT1] \S 4).
We define
$$\bar S_{\bar p}(\varepsilon):=\{{\bar p}+\beta|\beta\in
A^{01}\End_0(E)\oplus A^0E, \bar
D^0_{\bar p}{\bar {D_{\bar p}^0}}^*(\beta)+{\bar {D_{\bar
p}^1}}^*(\bar\sigma({\bar p}+\beta))=0,
||\beta||<\varepsilon\}.$$
The same arguments as in the proof of Proposition 3.9 show that for
sufficiently small
$\varepsilon>0$,
$\bar S_{\bar p}(\varepsilon)$ is a submanifold of $\bar{\cal A}$, whose
tangent space in $\bar p$
coincides with the first harmonic space $\bar\H^1_{\bar p}$ of the elliptic
complex  $(\bar{\cal
D}_{\bar p})$. Therefore,  we get a local finite dimensional model
$Z(\bar\sigma|_{\bar
S_{\bar p}(\varepsilon)}) $ for the intersection $Z(\bar\sigma)\cap \bar
N_{\bar p}(\varepsilon)$
of the integrable locus with the standard slice
$$ \bar N_{\bar p}(\varepsilon):=\{{\bar p} +\beta|\beta\in
A^{01}(\End_0(E))\oplus A^0(E),  {\bar{ D_{\bar p}^0}}^*(\beta)=0,\
||\beta||<\varepsilon\}$$
 through $\bar p$. The restriction
$$\bar\pi|_{Z(\bar\sigma|_{\bar S_{\bar p}(\varepsilon)})} :
Z(\bar\sigma|_{\bar
S_{\bar p}(\varepsilon)})\map {\cal M}^{s}_X(E,{\cal L})$$
of the orbit map is  \'etale if $[{\bar p}]\in {\cal M}^{ss}_X(E,{\cal
L})$, and induces an open
injection
$$\qmod{ Z(\bar\sigma|_{\bar S_{\bar p}(\varepsilon)})}{ \Z/r  }\map {\cal
M}^{s}_X(E,{\cal L})$$
if $[{\bar p}]\in {\cal M}^{s}_X(E,{\cal L})\setminus {\cal
M}^{ss}_X(E,{\cal L})$. We define the
orbifold structure of ${\cal M}^{s}_X(E,{\cal L})$ by taking the maps
$\bar\pi|_{Z(\bar\sigma|_{\bar S_{\bar p}(\varepsilon)})}$ as orbifold-charts.
\qed

Our next purpose is to compare the two types of moduli spaces constructed
in this paper.
Let $(X,g)$ be a K\"ahler surface endowed with the canonical
$Spin^c$-structure $\c $. Let
$E$ be a $U(2)$ bundle with $\det E=K_X$, and denote by ${\cal M}^*(E)$ the
moduli space
of irreducible quaternionic monopoles associated to the
$Spin^h(4)$-structure defined by
$(\c,E^{\vee})$ (Lemma 1.3)

It follows from Proposition 4.1     that ${\cal M}^*(E)$ has a
decomposition
$${\cal M}^*(E)= {\cal M}^*(E)_{\alpha=0}\cup  {\cal M}^*(E)_{\varphi=0}\ ,$$
where ${\cal M}^*(E)_{\alpha=0}$ ( ${\cal M}^*(E)_{\varphi=0}$ ) is the
Zariski closed
subspace of ${\cal M}^*(E)$ cut out by the equation $\alpha=0$ (
${\varphi=0}$ ). The
intersection
$${\cal M}^*(E)_{\alpha=0}\cap  {\cal M}^*(E)_{\varphi=0}$$
 is the Donaldson instanton  moduli space ${\cal D}^*$ of irreducible
projectively anti-self-dual
connections in $E$,  inducing the Chern connection in ${\cal K}_X$.

\begin{pr} The affine isomorphism ${\cal
A}\ni(C,\varphi)\longmapsto(\bar\partial_C,\varphi)\in\bar{\cal A}$ induces
a natural real
analytic open embedding
$$J:{\cal M}^*(E)_{\alpha=0}\hookrightarrow {\cal M}^s(E,{\cal K}_X)$$
whose image is the suborbifold of stable pairs of type $(E,{\cal K}_X)$.
\end{pr}
\pf Standard arguments (cf. [OT1]) show that $J$ is an \'etale map which
induces natural
identifications of the local models.

A point $[(\bar\delta,\varphi)]$  lies in the image of $J$ iff the
$SL(E)$-orbit of
$(\bar\delta,\varphi)$ intersects the zero locus of the map
$$m:\bar{\cal A}\map A^0(su(E)), \ \ (\bar\partial_C,\varphi)\longmapsto
\Lambda_g
F_C^0-\frac{1}{2}(\varphi\bar\varphi)_0 \ .$$
Let $({\cal E},\varphi)$ be the holomorphic pair of type $(E,{\cal K}_X)$
defined by
$(\bar\delta,\varphi)$. We can reformulate the condition above in the
following way:
$[({\cal E},\varphi)]$ lies in the image of $J$ iff there exists a
Hermitian metric $h$ in
${\cal E}$ inducing the K\"ahler metric in ${\cal K}_X=\det{\cal E}$ which
satisfies the
projective vortex equation $(V)$. But we know already that this holds iff
$({\cal
E},\varphi)$ is stable. Moreover, the unicity of the solution of the
projective vortex
equation is equivalent to the injectivity of $J$.
\qed

Using  the remark after Proposition 4.1, we can now state the  main result
of this
paper:
\begin{th}

Let $(X,g)$ be a K\"ahler surface with canonical bundle ${\cal K}_X$, and
let $E$ be a
$U(2)$-bundle with $\det E=K_X$. Consider the
$Spin^h$-structure associated with the canonical $Spin^c(4)$-structure and the
$U(2)$-bundle $E^{\vee}$.  The corresponding moduli space of irreducible
quaternionic
monopoles  is a union of two Zariski closed subspaces. Each of these
subspaces is
naturally isomorphic as a real analytic orbifold to the moduli space of
stable pairs of type
$(E,{\cal K}_X)$. There exists a real analytic involution on the
quaternionic moduli space
which interchanges these two closed subspaces. The fixed point set of this
involution is
 the Donaldson moduli space of instantons in $E$ with fixed determinant, modulo
the gauge group $SU(E)$. The closure of the complement of the Donaldson
moduli space
intersects the moduli space of instantons in the Brill-Noether locus.

The union ${\cal M}^{SW}$ of all rank 1-Seiberg-Witten moduli spaces
associated with
splittings   $E=E' \oplus E''$   corresponds to the subspace of stable
pairs of type ii).
\end{th}

\newpage

\centerline{\large{\bf References}}
\vspace{10 mm}
\parindent 0 cm

[B1] Bradlow, S. B.: {\it Vortices in holomorphic line bundles over closed
K\"ahler manifolds}, Comm. Math. Phys. 135, 1-17 (1990)

[B2] Bradlow, S. B.: {\it Special metrics and stability for holomorphic
bundles with global sections}, J. Diff. Geom. 33, 169-214 (1991)

[D] Donaldson, S.: {\it Anti-self-dual Yang-Mills connections over complex
algebraic surfaces and stable vector bundles}, Proc. London Math. Soc. 3,
1-26 (1985)

[DK] Donaldson, S.; Kronheimer, P.B.: {\it The Geometry of four-manifolds},
Oxford Science Publications (1990)

[FU] Freed D. S. ;  Uhlenbeck, K.:
{\it Instantons and Four-Manifolds.}
Springer-Verlag 1984.

[HH] Hirzebruch, F.; Hopf H.: {\it Felder von Fl\"achenelementen in
4-dimensionalen 4-Mannigfaltigkeiten}, Math. Ann. 136 (1958)

[H] Hitchin, N.: {\it  Harmonic spinors}, Adv. in Math. 14, 1-55 (1974)

[JPW] Jost, J.; Peng, X.;  Wang, G. :{\it Variational aspects of the
Seiberg-Witten functional},
Preprint, dg-ga/9504003, April (1995)

[K] Kobayashi, S.: {\it Differential geometry of complex vector bundles},
Princeton University Press (1987)

[KM] Kronheimer, P.; Mrowka, T.: {\it The genus of embedded surfaces in the
projective plane}, Preprint (1994)

[LM] Labastida, J. M. F.;  Marino M.: {\it Non-abelian monopoles on four
manifolds}, Preprint,
Departamento de Fisica de Particulas,   Santiago de Compostela, April
 (1995)

[LT] L\"ubke, M.; Teleman, A.: {\it The Kobayashi-Hitchin correspondence},
 World Scientific Publishing Co,  to appear.

[M]  Miyajima K.: {\it Kuranishi families of
vector bundles and algebraic description of
the moduli space of Einstein-Hermitian
connections},   Publ. R.I.M.S. Kyoto Univ.  25,
  301-320 (1989)

[OSS] Okonek, Ch.; Schneider, M.; Spindler, H: {\it Vector bundles  on complex
projective spaces}, Progress in Math. 3, Birkh\"auser, Boston (1980)

[OT1]   Okonek, Ch.; Teleman A.: {\it The Coupled Seiberg-Witten Equations,
Vortices, and Moduli Spaces of Stable Pairs}, Preprint, Z\"urich, Jan. 13-th,
(1995)

[OT2] Okonek, Ch.; Teleman A.: {\it Les invariants de Seiberg-Witten et la
conjecture de Van De  Ven}, to appear in Comptes Rendus

[OT3] Okonek, Ch.; Teleman A.: {\it Seiberg-Witten invariants and rationality
of complex surfaces}, Preprint, Z\"urich, March (1995)

[W] Witten, E.: {\it Monopoles and four-manifolds}, Mathematical  Research
Letters 1,  769-796  (1994)
\vspace{0.4cm}\\
Authors addresses:\\
\\
Mathematisches Institut, Universit\"at Z\"urich,\\
Winterthurerstrasse 190, CH-8057 Z\"urich\\
e-mail: okonek@math.unizh.ch  \ ; \
 teleman@math.unizh.ch

\end{document}